\documentclass[11pt]{article}

\usepackage{amsmath, amsfonts, amssymb, amsthm}

\usepackage{setspace}

\DeclareMathSymbol{\leqslant}{\mathalpha}{AMSa}{"36} 

\DeclareMathSymbol{\geqslant}{\mathalpha}{AMSa}{"3E} 

\DeclareMathSymbol{\eset}{\mathalpha}{AMSb}{"3F}     

\newcommand{\Om}{\Omega}

\newcommand{\IR}{\mathbb{R}}

\newcommand{\IE}{\mathbb{E}}

\newcommand{\ban}{\begin{align}}
\newcommand{\ean}{\end{align}}

\newcommand{\ba}{\begin{align*}}
\newcommand{\ea}{\end{align*}}

\newcommand{\be}{\begin{eqnarray*}}
\newcommand{\ee}{\end{eqnarray*}}

\newcommand{\ben}{\begin{eqnarray}}
\newcommand{\een}{\end{eqnarray}}

\theoremstyle{plain}
\newtheorem{theo}{Theorem}[section]

\newtheorem{propo}[theo]{Proposition}
\newtheorem{corollary}[theo]{Corollary}

\newtheorem{remark}[theo]{Remark}
\theoremstyle{definition}

\begin{document}

\vglue20pt \centerline{\huge\bf Uniform shrinking and expansion }

\medskip

\centerline{\huge\bf  }

\medskip

\centerline{\huge\bf  under isotropic  Brownian flows }
 \medskip


\bigskip

\bigskip

\centerline{by}

\bigskip

\medskip

\centerline{{\Large Peter Baxendale\footnotemark[1] and  Georgi  Dimitroff\footnotemark[2]}}

\footnotetext[1]{Department of Mathematics, University of Southern
California, 3620 S. Vermont Avenue, Los Angeles, CA 90089-2532.
Supported in part by NSF Grant DMS-05-04853.}
\footnotetext[2]{Fraunhofer ITWM, Fraunhofer-Platz 1, D-67663 Kaiserslautern}

\bigskip

\bigskip

{\leftskip=1truecm

\rightskip=1truecm

\baselineskip=15pt

\small

\noindent{\slshape\bfseries Summary.} We study some finite time
transport properties of isotropic Brownian flows. Under a certain
nondegeneracy condition on the potential spectral measure, we
prove that uniform shrinking or expansion of balls under the flow
over some bounded time interval can happen with positive
probability.  We also provide a control theorem for isotropic
Brownian flows with drift. Finally, we apply the above results to
show that under the nondegeneracy condition the length of a
rectifiable curve evolving in an isotropic Brownian flow with
strictly negative top Lyapunov exponent converges to zero as $t\to
\infty$ with positive probability.
 \bigskip


\bigskip


}


\newcommand{\oq}{{\langle}}
\newcommand{\ve}{{\varepsilon}}
\newcommand{\cq}{{\rangle}_t}
\newcommand{\dd}{{\mathrm{d}}}
\newcommand{\Dd}{{\mathrm{D}}}
\newcommand{\cd}{{\cdot}}
\newcommand{\nix}{{\varnothing}}
\newcommand{\N}{{\mathbb N}}
\newcommand{\Z}{{\mathbb Z}}
\newcommand{\R}{{\mathbb R}}
\newcommand{\Q}{{\mathbb Q}}
\newcommand{\E}{{\mathbb E}}
\renewcommand{\P}{{\mathbb P}}
\newcommand{\F}{{\cal F}}
\newcommand{\C}{{\mathbb C}}
\newcommand{\K}{{\mathbb K}}
\newcommand{\B}{{\cal B}}
\newcommand{\G}{{\cal G}}
\newcommand{\D}{{\cal D}}
\newcommand{\Zz}{{\cal Z}}
\newcommand{\Ll}{{\cal L}}
\newcommand{\A}{{\cal A}}
\newcommand{\Po}{{\cal P}}
\newcommand{\Sy}{{\cal S}}
\newcommand{\cZ}{{\cal Z}}
\newcommand{\M}{{\cal M}}
\newcommand{\Nn}{{\cal N}}
\newcommand{\p}{{\mathbf P}}
\newcommand{\X}{{\mathbb X}}
\newcommand{\interior}[1]{{\mathaccent 23 #1}}
\newcommand{\bem}{\begin{em}}
\newcommand{\eem}{\end{em}}
\def\eins{{\mathchoice {1\mskip-4mu\mathrm l}
{1\mskip-4mu\mathrm l}
{1\mskip-4.5mu\mathrm l} {1\mskip-5mu\mathrm l}}}
\newcommand{\olim}[1]{\begin{array}{c} ~\\[-1.4ex] \overline{\lim} \\[-1.35ex]
                       {\scriptstyle #1}\end{array}}
\newcommand{\plim}[2]{\begin{array}{c} #1\\[-1.175ex] \longrightarrow \\[-1.2ex]
                       {\scriptstyle #2}\end{array}}

\newcommand{\ulim}[1]{\begin{array}{c} ~\\[-1.175ex] \underline{\lim} \\[-1.2ex]
                       {\scriptstyle #1}\end{array}}

\renewcommand{\theequation}{\thesection.\arabic{equation}}

\section{Introduction}\label{sec-intro}
Stochastic flows generalize the notion of stochastic process and it  has been suggested that they are a natural probabilistic model for the evolution of  passive tracers within a turbulent fluid.
  The one point motions (trajectories of individual particles) are diffusions and the motions of adjacent points are correlated
and form a stochastic flow of homeomorphisms. The flow is
generated by a random  field of continuous semimartingales
$F(t,x)$ in the sense that loosely speaking for a short  time
increment $t\to t+\Delta t$ the passive tracer starting in $x\in
\R^d$ at time $t$ moves along the increment of the field $F$, i.e.
we have
\begin{align*}
\phi_{t,t+\Delta t}(x)-x \approx F(t+\Delta t, x)-F(t,x).
\end{align*}
Infinitesimally, the above dynamics is achieved by the family of
solutions (for different initial conditions) of the stochastic
differential equation of Kunita type
\begin{align*}
\phi_{s,t}(x)=x+\int\limits_s^t F(\dd u, \phi_{s,u}(x)).
\end{align*}
We consider isotropic Brownian flows (IBF) on $\R^d$ for $d \ge
2$.  These are a special class of stochastic flows characterized
(modulo regularity conditions) by their spatial translation and
rotational invariance, temporal homogeneity and independence of
their increments. They have been extensively studied by many
authors, including Baxendale and Harris \cite{BaH86}, Le Jan
\cite{LeJ85}, Le Jan and Darling \cite{DaLe88}, and Cranston,
Scheutzow and Steinsaltz \cite{CSS99}. In Section \ref{sec-IBF} we
will give a short introduction to IBF with emphasis on properties
and facts needed for our purposes later on.

Often when speaking about macroscopic properties of stochastic
flows one illustrates the action of the flow by the evolution of
an oil spill on the surface of an ocean, where the oil spill is a
passive tracer following the turbulence of the ocean's surface. In
Scheutzow and Steinsaltz \cite{SS02} (see also \cite{CSS99} and
\cite{CSS00}) it has been shown that the diameter of the oil spill
grows linearly in time as $t\to \infty$,  almost surely if the top
Lyapunov exponent $\lambda$ is positive, and with strictly
positive probability if $\lambda \le 0$.  In this paper, instead
of looking at asymptotic behavior of the IBF as $t \to \infty$, we
will study the behavior at a fixed time $T$ or, more generally,
during a fixed time interval $[T_1,T_2]$.   We consider the
following question: is there a time interval $[T_1,T_2]$ such that
a circular oil spill is uniformly squeezed with positive
probability by the action of the ocean during the interval
$[T_1,T_2]$?  More precisely, we will be discussing the following
question: Does an isotropic Brownian flow squeeze a ball of radius
$R$ into a ball of radius $r < R$ with positive probability during
some fixed time interval $[T_1,T_2]$?

Since the radial component of the motion of every point on the
boundary of the ball of radius $R$ is a time changed (scalar)
Brownian motion, the motion of every boundary point will almost
surely cross the boundary infinitely many times during any time
interval $(0,t]$. Therefore the ball cannot  be mapped into
itself throughout any interval $(0,t]$, and we will restrict
attention to intervals $[T_1,T_2]$ with $T_1 > 0$.

Obviously the answer  to our question is ``No'' if the flow is
volume-preserving.  However it turns out that apart from the
volume-preserving case, the answer is usually ``Yes'', subject to
a particular non-degeneracy condition on the potential spectral
measure $M_P$ associated with the flow.  A brief survey of
isotropic Brownian flows, including the definition of $M_P$, is
given in Section \ref{sec-IBF}.  The non-degeneracy condition
(condition $\mathbf{(C)_\rho}$) is given in Section \ref{sec-main}
along with the main results (Theorems \ref{mt_thm1} and
\ref{mt_thm2} and Corollary \ref{cor1}) on the squeezing of balls.
Moreover a simple adaptation of the proofs shows that, under
exactly the same condition, the flow can also expand balls (see
Remark \ref{rem-exp}).

The proofs remain valid if we add a deterministic drift.  For the
sake of generality we do this and think of the drift part as the
deterministic current in the ocean and the random isotropic
Brownian part as the unpredictable turbulent movement on the
surface.  The drift is assumed to be time homogeneous for sake of
simplicity.

Section \ref{sec-rkhs} contains material on the reproducing kernel
Hilbert space associated with the Brownian field $M(t,x)$, and
Section \ref{sec-control} contains a control theorem for
stochastic flows.  Both of these sections may be of independent
interest.  The proofs of the main theorems are given in Section
\ref{sec-pf}.  Finally, Section \ref{sec-length} considers how the
length of a curve evolves under an isotropic Brownian flow, and
extends a result of Baxendale and Harris \cite{BaH86}.

\section{Isotropic Brownian flows} \label{sec-IBF}

Here we  provide a short  introduction to IBF following mainly
Baxendale and Harris \cite{BaH86}, Le Jan \cite{LeJ85} and Yaglom
\cite{Yag57}.

A \bem (forward) stochastic flow of homeomorphisms \eem on $\R^d$
is a family of random homeomorphisms  $\left\{\phi_{s,t} : 0\le
s\le t< \infty\right\}$ of $\R^d$ into itself, such that almost
surely $\phi_{u,t}\circ \phi_{s,u}=\phi_{s,t}$ for $s\le u\le t$
and $\phi_{t,t}=\text{Id}_{\R^d}$.  If the increments $\phi_{s,t}$
on disjoint intervals are independent and time homogeneous then
the flow is said to be a \bem Brownian flow. \eem

According to Kunita \cite[Theorem 4.2.8]{Ku90}, under suitable
regularity conditions, Brownian flows of homeomorphisms can be
realized as solutions of Kunita-type  SDEs
\begin{align}\label{pre_ibf1}
\phi_{s,t}(x)=x+\int \limits_s^tM(\dd u,\phi_{s,u}(x)) + \int_0^t
v(\phi_{s,u})\dd u.
\end{align}
Here $v:\IR^d \to \IR^d$ is a vector field and $M:\IR_+\times
\IR^d\times \Om\to \IR^d$ is a mean-zero Gaussian random field.
$M$ is called the generating \bem Brownian field \eem and its
distribution is determined by the covariances
$$
\IE \left[\langle M(t,x),\xi\rangle \langle M(s,y),\eta\rangle
\right] =(s\wedge t) \,\langle b(x,y)\xi,\eta\rangle,   \quad \quad
\xi,\eta \in \IR^d,
$$
for some covariance tensor $b: \IR^d \times \IR^d \to L(\IR^d)$.
The function $b$ is positive semi-definite:  for all $n \ge 1$,
all $x^{(1)},x^{(2)}, \ldots,x^{(n)} \in \IR^d$ and all
$\xi^{(1)}, \xi^{(2)}, \ldots, \xi^{(n)} \in \IR^d$ we have
  \begin{equation} \label{bpos}
     \sum_{k,\ell = 1}^n \left\langle
     b(x^{(k)},x^{(\ell)})\xi^{(k)},\xi^{ (\ell)} \right\rangle \ge 0.
   \end{equation}
The law of the stochastic flow $\{\phi_{s,t}: 0 \le s \le t <
\infty\}$ is determined by the functions $b(x,y)$ and $v(x)$. The
differentiability properties of the mappings $\phi_{s,t}$ depend on
the differentiability of the functions $b(x,y)$ and $v(x)$.  \\

An \emph{isotropic Brownian flow}  (IBF) on $\IR^d$ is a Brownian
flow of diffeomorphisms of $\IR^d$ for which the distribution of
each $\phi_{s,t}$ is invariant under rigid transformations of
$\IR^d$.  IBF have been extensively studied by Baxendale and
Harris \cite{BaH86} and Le Jan \cite{LeJ85}. The invariance in
distribution of $\phi_{s,t}$ under rigid motions implies the
invariance in distribution of the generating Brownian field
$M(t,x)$; we say that $M$ is an {\em isotropic Brownian field}.
The invariance under translations implies that $b(x,y) = b(x-y,0)
\equiv b(x-y)$ and then the invariance under rotations and
reflections implies that
 \begin{equation} \label{bisot}
   b(x)=O^Tb(Ox)O
    \end{equation}
for all $O$ in the orthogonal group $\mathcal{O}^d$. Moreover for
IBF we have $v(x) \equiv 0$.  In this paper we will assume that
$b(x)$ is $C^4$.  Then $b(x-y)$ will be $C^{2,2}$ as a function of
$x$ and $y$, and so the resulting IBF will consist of $C^1$
diffeomorphisms, see Kunita \cite[Theorem 4.6.5]{Ku90}. Moreover,
the isotropy property (\ref{bisot}) implies that $b(0) = cI$ for
some constant $c$.  At the cost of rescaling time by a constant
factor we can and will assume that $b(0) = I$.  In order to avoid
the trivial case where the flow consists of translations, we
assume also that $b(x) \not\equiv I$.
\\

According to Yaglom \cite[Section 4]{Yag57} (and as described in
\cite{BaH86}) a positive semi-definite tensor $b(x)$ with the
above properties can be written in the form
\begin{equation}\label{pre_ibf2}
b_{ij}(x)=
 \big(B_L(|x|)-B_N(|x|)\big)\frac{x_ix_j}{|x|^2}+\delta_{ij}B_N(|x|)
\end{equation}
for $x\ne 0$, where $B_L$ and $B_N$ are the so-called longitudinal
and transverse covariance  functions defined by
$$
B_L(r)=b_{ii}(re_i),\quad  \quad  B_N(r)=b_{ii}(re_j)
 $$
for $r \ge 0$ and any $i \neq j$.  As usual $e_i$ denotes the
$i$-th standard basis vector in $\R^d$. $B_L$ and $B_N$ are
bounded $C^4$ functions with bounded derivatives. Further, the
isotropic covariance tensor $b$ can be decomposed
  \begin{align}\label{pre_ibf_decom}
b(x)=\mu_0 I+\mu_1b_P(x)+\mu_2b_S(x)
\end{align}
with $\mu_i \ge 0$ and $\mu_0+\mu_1+\mu_2 = 1$, where $b_P(x)$ is
the covariance tensor for an isotropic Brownian field consisting
of gradient vector fields, and $b_S(x)$ is the covariance tensor
for an isotropic Brownian field consisting of divergence-free
vector fields.  The labels $P$ and $S$ stand for ``potential'' and
``solenoidal'' respectively.  With the normalizing conditions
$b^P(0) = b^S(0) = I$ and $\lim_{|x| \to \infty} b_P(x) =
\lim_{|x| \to \infty} b_S(x) = 0$ the decomposition is unique
(except trivially when $\mu_1 = 0$ or $\mu_2 = 0$).  The functions
$b_P$ and $b_S$ can each be written in the form (\ref{pre_ibf2})
and the corresponding longitudinal and transverse covariance
functions $B_{PL}$, $B_{PN}$, $B_{SL}$ and $B_{SN}$ are uniquely
determined by two finite spectral measures $M_P$ and $M_S$ on the
positive real line $(0,\infty)$ through the expressions
\begin{align}\label{pre_ibf_spect}
\begin{split}
&B_{PL}(s)=2^{\frac{d-2}{2}}\Gamma\left(\frac{d}{2}\right)\int\limits_{(0,\infty)}\left[
  \frac{\mathcal{J}_{\frac{d}{2}}(sr)}{(sr)^{\frac{d}{2}}}-
  \frac{\mathcal{J}_{\frac{d+2}{2}}(sr)}{(sr)^{\frac{d-2}{2}}}\right]M_P(\dd r) \,,\\
&B_{PN}(s)=2^{\frac{d-2}{2}}\Gamma\left(\frac{d}{2}\right)\int\limits_{(0,\infty)}
\frac{\mathcal{J}_{\frac{d}{2}}(sr)}{(sr)^{\frac{d}{2}}}M_P(\dd r) \,,\\
&B_{SL}(s)=2^{\frac{d-2}{2}}\Gamma\left(\frac{d}{2}\right)(d-1)\int\limits_{(0,\infty)}
  \frac{\mathcal{J}_{\frac{d}{2}}(sr)}{(sr)^{\frac{d}{2}}}M_S(\dd r) \,,\\
&B_{SN}(s)=2^{\frac{d-2}{2}}\Gamma\left(\frac{d}{2}\right)\int\limits_{(0,\infty)}\left[
  \frac{\mathcal{J}_{\frac{d-2}{2}}(sr)}{(sr)^{\frac{d-2}{2}}}-
  \frac{\mathcal{J}_{\frac{d}{2}}(sr)}{(sr)^{\frac{d}{2}}}\right]M_S(\dd r),
\end{split}
\end{align}
where $\mathcal{J}_{\nu}$ denotes the Bessel function (of the
first kind) of order $\nu$. $M_P$ and $M_S$ are the {\em potential
spectral measure} and the {\em solenoidal spectral measure}
respectively for the isotropic covariance function $b$. The
normalization $b_P(0) = b_S(0) = I$ gives $M_P((0,+\infty))=d$ and
$M_S((0,+\infty))=\frac{d}{d-1}$, and the assumption that $b$ is
$C^4$ implies that $M_P$ and $M_S$ both have finite $4^{th}$
moments. For future reference, write
  \begin{equation}
\label{betaL}
  \beta_L := -\frac{\partial^2b_{ii}}{\partial x_i^2}(0) = -B_L''(0) = \frac{3\mu_1}{d(d+2)}\int s^2 M_P(\dd s) +
  \frac{(d-1)\mu_2}{d(d+2)}\int s^2M_S(\dd s) > 0
    \end{equation}
 and
  \begin{equation} \label{betaN}
  \beta_N := -\frac{\partial^2b_{ii}}{\partial x_j^2}(0) =-B_N''(0) = \frac{\mu_1}{d(d+2)}\int s^2 M_P(\dd s) +
  \frac{(d+1)\mu_2}{d(d+2)}\int s^2M_S(\dd s) > 0
    \end{equation}
for $j \neq i$.  Conversely, to every pair of suitably normalized
$M_P$ and $M_s$ with finite $4^{th}$ moments and non-negative
constants $\mu_1+\mu_2 \le 1$ we can construct a $C^4$ isotropic
covariance tensor $b$ and consequently an isotropic Brownian field
and (via a Kunita-type SDE) also an isotropic Brownian flow.

\section{Uniform shrinking and expansion} \label{sec-main}

As mentioned in the introduction we will give the main result in
the more general setting of  IBF with drift rather than for pure
IBF. Let $M(t,x,\omega)\colon \R_+\times \R^d\times \Omega \to
\R^d$ be an isotropic Brownian generating field with a $C^4$
covariance tensor $b$ satisfying $b(0) = I$ and $b(x) \not\equiv
I$. Let $v(x):\R^d \to \R^d $ be a deterministic $C_{b}^{2}$
function (derivatives of order up to $2$ exist and are bounded and
continuous). Consider the semimartingale field
$$F(t,x)=M(t,x)+tv(x)$$
and the stochastic flow $\phi$ generated via the SDE:
$$\phi_{s,t}(x)=x+\int\limits_s^t F(\dd u,\phi_{s,u}(x))=x+\int\limits_s^t M(\dd u,\phi_{s,u}(x))+\int\limits_s^tv(\phi_{s,u}(x))\dd u.$$
We will call $\phi$ an \bem isotropic Brownian flow with drift.
\eem Of course, in case $v\equiv 0$,  $\phi$ is simply an IBF.
Before we proceed with the main result we shall state the
condition on the isotropic covariance tensor $b(x)$ which will be
needed below.  Suppose $\rho > 0$

\vspace{2ex}

\noindent {\bf Condition $\mathbf{(C)_\rho}$:}  $\mu_1
> 0$ in the decomposition (\ref{pre_ibf_decom}) and the
potential spectral measure $M_P$ is not supported on the set of
all zeros of the mapping $s \mapsto \mathcal{J}_\frac{d}{2}(\rho
s)$, i.e.
supp$(M_P)\not \subset \{s \ge 0: \mathcal{J}_\frac{d}{2}(\rho s) = 0\}$. \\

\begin{remark}{\rm  (i) The condition $\mu_1 > 0$ says
that the isotropic Brownian field $M(t,x)$ is not divergence-free.
From (\ref{betaL}) and (\ref{betaN}) we have $(d+1)\beta_L -
(d-1)\beta_N \ge 0$, with equality if and only $\mu_1 = 0$. Recall that the flow 
is divergence-free (incompressible) if and only if $(d+1)\beta_L -
(d-1)\beta_N = 0$.

(ii) If $\mu_1 > 0$ and the potential spectral measure $M_P$ has a
density then ${\bf (C)}_\rho$ is satisfied for all $\rho > 0$.
Using Fourier transform theory, the covariance function $b$ can be
written
   $$
   b(x) = \int_{\IR^d} e^{i\langle x,\lambda \rangle} F(\dd \lambda)
   $$
where $F$ is a $L(\IR^d)$ valued measure with the property that
$F(A)$ is non-negative definite for every Borel subset $A \subset
\IR^d$.  The isotropy property (\ref{bisot}) of $b$ implies that
$F$ can be written with respect to polar coordinates $r =
|\lambda|$ and $\theta = \lambda/|\lambda|$ in the form
   $$
   F(\dd \theta, \dd r) = \mu_1\theta_j \theta_k \sigma(\dd \theta)M_P(\dd r) +
   \mu_2 (\delta_{jk}-\theta_j\theta_k) \sigma(\dd \theta) M_S(\dd r)
   $$
for $r > 0$, where $\sigma$ denotes the uniform probability
measure on the unit sphere $\mathbb{S}^{d-1}$. (This is the
argument used by Yaglom \cite{Yag57} to derive the formulas
(\ref{pre_ibf_spect}).)  Therefore the integrability condition
  $$
  \int_{\IR(d)} \|b(x)\|\dd x < \infty
  $$
implies the existence of a density for $F$, which in turn implies
the existence of a density for $M_P$ and hence condition
$\mathbf{(C)_\rho}$ for all $\rho > 0$, so long as $\mu_1
> 0$.}
\end{remark}

The first main result deals with the squeezing, or contraction, of
a ball of some fixed radius $R$.  Since the distribution of the
IBF is invariant under translations, it is enough to consider the
ball ${\text B}(0,R)$ of radius $R$ centered at 0.

\begin{theo}\label{mt_thm1} Let $\big\{\phi_{s,t}(x) \colon 0 \le
s \le t < \infty,\, x\in \R^d\big\}$ be an isotropic Brownian flow
with drift with generating field $\big\{F(t,x) \colon t\ge
0,\,x\in \R^d\big\}$. For $R>0$, assume that the covariance tensor
$b$ satisfies condition $\mathbf{(C)}_R$. Then there exists
$\delta
>0$  such that
\begin{equation} \label{Rcont}
\P\big(\phi_t\left(\textnormal{B}(0,R+\delta)\right)\subset
\textnormal{B}(0,R-\delta) \mbox{ \rm{for all} } t \in[T_1,T_2]
\big)
> 0
\end{equation}
for all $0 < T_1 < T_2$.
\end{theo}

The second main result extends this result to the squeezing of the
closed ball $\overline{\text B}(0,R)$ inside the open ball ${\text
B}(0,r)$ for arbitrary $R$ and $r$.  Clearly it is enough to
consider the case where $r < R$.

\begin{theo} \label{mt_thm2} Let $\big\{\phi_{s,t}(x) \colon 0 \le
s \le t < \infty,\, x\in \R^d\big\}$ be an isotropic Brownian flow
with drift with generating field $\big\{F(t,x) \colon t\ge
0,\,x\in \R^d\big\}$. Suppose $0 < r < R$.  Assume that the
covariance tensor $b$ satisfies condition $\mathbf{(C)_\rho}$ for
all $\rho \in [r,R]$. Then
\begin{equation} \label{rRcont}
\P \big(\phi_{t}\left(
\overline{\textnormal{B}}(0,R)\right)\subset \mathrm{B}(0,r)
\mbox{ \rm{for all}  } t \in [T_1,T_2] \big)>0\,
\end{equation}
for all $0 < T_1 < T_2$.
\end{theo}

\begin{corollary}\label{cor1}
Let $\big\{\phi_{s,t}(x) \colon 0 \le s \le t < \infty,\, x\in
\R^d\big\}$ be an isotropic Brownian flow with drift with
generating field $\big\{F(t,x)  \colon t\ge 0,\,x\in \R^d\big\}$.
Assume that the covariance tensor satisfies condition
$\mathbf{(C)_\rho}$ for all $\rho
> 0$. Then for arbitrary $R,r>0$ the inequality (\ref{rRcont}) holds for all $0 < T_1 < T_2$.
\end{corollary}

Clearly Corollary \ref{cor1} is an immediate consequence of
Theorem \ref{mt_thm2}.  The proofs of Theorems \ref{mt_thm1} and
\ref{mt_thm2} will be given in Section \ref{sec-pf}.

\begin{remark}\label{rem-exp} {\rm  There are versions of these three
results asserting expansion instead of expansion.  Under exactly
the same conditions, the equations (\ref{Rcont}) and
(\ref{rRcont}) can be replaced by
   \begin{equation} \label{Rexp}
\P\big(\phi_t\left(\textnormal{B}(0,R-\delta)\right)\supset
\textnormal{B}(0,R+\delta) \mbox{ \rm{for all} } t \in[T_1,T_2]
\big)
> 0
\end{equation}
and
 \begin{equation} \label{rRexp}
\P \big(\phi_{t}\left( \textnormal{B}(0,r)\right)\supset
\overline{\mathrm{B}}(0,R) \mbox{ \rm{for all}  } t \in [T_1,T_2]
\big)>0\,
\end{equation}
respectively.}
\end{remark}

\section{RKHS and the generating Brownian field} \label{sec-rkhs}

Recall the covariance function $b(x,y) \in L(\IR^d)$ of the
generating Brownian field $M(t,x)$ is given by
   $$
   \langle b(x,y)\xi,\eta\rangle = \IE[\langle M(1,x),\xi\rangle \langle M(1,y),\eta \rangle
   ]\quad \quad \xi,\eta \in \IR^d.
   $$
Assume that $b$ is continuous in both variables and bounded.
Associated to the positive semi-definite function $b$ is a Hilbert
space $H$ consisting of vector fields on $\R^d$, that is $H
\subset C(\IR^d : \IR^d)$. The space $H$ is characterized by the
properties

(i) for each $x\in\mathbb{R}^d$ and $\xi \in \mathbb{R}^d$ the
vector field $b_{x,\xi}: y \to b(x,y)\xi$ is an element of $H$;
and

(ii)  $\langle f,b_{x,\xi}\rangle_{H} = \langle f(x),\xi \rangle$
for all $f \in H$ and $x \in \mathbb{R}^d$, $\xi \in\mathbb{R}^d$.

\noindent The property (ii) is called the reproducing property,
and says that taking an inner product with $b_{x,\xi}$ acts like
an evaluation map on $H$.  In this setting the function $b$ is
called a reproducing kernel and $H$ is the associated reproducing
kernel Hilbert space (RKHS).   For details of the theory of
reproducing kernel Hilbert spaces for vector valued functions see
Baxendale \cite[Section 5]{Bax76}.

The reproducing property (ii) shows how the Hilbert space $H$
determines the functions $b_{x,\xi}$ for all $x \in \IR^d$ and
$\xi \in\IR^d$.  These in turn determine the function $b$ and
hence the distribution of the generating Brownian field $M(t,x)$.
In particular $H$ is the Cameron Martin space for the distribution
$\gamma = \P \circ M(1,\cdot)^{-1}$, see Bogachev \cite{Bog98}.

The main result in this section is Theorem \ref{thm-V}, which
shows that the condition $\mathbf{(C)_\rho}$ is equivalent to the
existence of a particular sort of vector field in the RKHS $H$.

\begin{propo} \label{prop-int1} Let $H$ be a reproducing kernel Hilbert space of functions
$\IR^d \to \IR^d$, and suppose $\mu$ is a measure on $H$ such that
    $$
    \int \|U\|_H \,\dd\mu(U) < \infty
    $$
Then there exists $D  \in H$ such that
\begin{equation} \label{Ddef}
   \int \langle V, U \rangle_H\,\dd\mu(U)  = \langle V, D \rangle_H
   \quad \quad \mbox{ for all }V \in H.
    \end{equation}
Moreover
  $$
  D(x)  = \int U(x) \,\dd\mu(U) \quad \quad \mbox{ for all }x \in {\R^d}
  $$
and
   $$
   \|D\|_H^2 = \int \int \langle U_1,U_2 \rangle_H \,\dd\mu(U_1)d\mu(U_2).
   $$
   \end{propo}

\noindent{\bf Proof.}  Consider the linear mapping $H \to
\mathbb{R}$ given by
   $$
   V \mapsto \int \langle V,U \rangle_H \,\dd\mu(U).
   $$
Since
   $$
  \int \left|\langle V ,U
   \rangle_H \right| \,\dd\mu(U) \le  \int \|V\|_H \|U\|_H \,\dd\mu(U)  =  \|V\|_H\int  \|U\|_H \,\dd\mu(U)
   $$
the mapping is well-defined and continuous and so there exists $D
\in H$ satisfying (\ref{Ddef}).  Taking $V = b_{x,\xi}$ we get
   $$
   \langle D(x),\xi \rangle = \langle D, b_{x,\xi}
   \rangle_H
    =\int \langle U, b_{x,\xi} \rangle_H\,\dd\mu(U)
     =  \int \langle U(x),\xi \rangle\,\dd\mu(U)= \left\langle \int U(x)\,\dd\mu(U),\xi
     \right\rangle
  $$
so that
    $$
    D(x) =  \int U(x)\,\dd\mu(U).
    $$
Finally
 \begin{eqnarray*}
  \|D\|^2_H = \langle D,D \rangle_H
   = \int \langle D,U_2 \rangle_H \,\dd\mu(U_2)
   & = &  \int \left( \int \langle U_1,U_2 \rangle_H\,\dd\mu(U_1)
   \right) \,\dd\mu(U_2)\\
   & = & \int \int\langle U_1,U_2 \rangle_H \,\dd\mu(U_1)\dd\mu(U_2)
  \end{eqnarray*}
and the proof is complete.  \hfill $\Box$

\vspace{2ex}

The result above does not use isotropy or even homogeneity.  For
the remainder of this subsection we restrict to the isotropic
situation.

\begin{theo}\label{thm-V}
Let $H$ be the reproducing kernel Hilbert space corresponding to
the covariance tensor $b(x)$ for an isotropic Brownian field
$M(t,x)$, and let $\rho
> 0$. There exists a vector field $V \in H$ which points inwards on
the boundary of the ball with radius $\rho$ and centered at the
origin, that is,
\begin{equation} \label{Vin}
\langle V(\rho \theta),\theta \rangle  < 0\,\,\,\,\,\textrm{
whenever  } \,\,\|\theta\| =1,
\end{equation}
if and only if the covariance tensor $b(x)$ satisfies condition
$\mathbf{(C)_\rho}$.
\end{theo}

\noindent \textbf{Proof.}  Denote by $\sigma$ the uniform
probability measure on the unit sphere $\mathbb{S}^{d-1}$, and let
$\mu$ be the probability measure which is the image of $\sigma$
under the mapping $\phi \mapsto b_{\rho\phi,\phi}$ of
$\mathbb{S}^{d-1}$ into $H$. Since $\phi \mapsto
\|b_{\rho\phi,\phi}\|_H = \sqrt{\langle b(0)\phi,\phi \rangle}$ is
bounded on $\mathbb{S}^{d-1}$ we can apply Proposition
\ref{prop-int1} to obtain a vector field $\widetilde{V} \in H$
given by
  \begin{equation} \label{D}
   \widetilde{V}(x) = \int_{\mathbb{S}^{d-1}} b_{\rho
   \phi,\phi}(x)\,\dd\sigma(\phi).
  \end{equation}
Using the isotropy of $b$ and the fact that $\sigma$ is invariant
under the action of $R \in {\cal SO}(d)$ we obtain for $\theta \in
\mathbb{S}^{s-1}$
  \begin{eqnarray*}
  \langle \widetilde{V}(\rho \theta),\theta \rangle
   & = & \int_{\mathbb{S}^{d-1}}\langle b_{\rho
   \phi,\phi}(\rho \theta), \theta \rangle)\,\dd\sigma(\phi)\\
    & = & \int_{\mathbb{S}^{d-1}}\langle b(\rho \phi - \rho \theta)\phi, \theta \rangle)\,\dd\sigma(\phi)\\
     & = & \int_{\mathbb{S}^{d-1}}\langle b(\rho R\phi - \rho R\theta)R\phi, R\theta \rangle)\,\dd\sigma(\phi)\\
      & = & \int_{\mathbb{S}^{d-1}}\langle b(\rho \phi - \rho R \theta)\phi, R\theta \rangle)\,\dd\sigma(\phi)\\
      & = & \langle \widetilde{V}(\rho R\theta),R\theta \rangle,
      \end{eqnarray*}
so that $\theta \mapsto \langle \widetilde{V}(\rho \theta),\theta
\rangle$ is constant.  Moreover, for all $U \in H$ we have
   $$
 \int_{\mathbb{S}^{d-1}} \langle U(\rho \phi),\phi \rangle
 \dd\sigma(\phi) = \int_{\mathbb{S}^{d-1}} \langle U,b_{\rho \phi,\phi}
 \rangle_H
 \dd\sigma(\phi) =\langle U,\widetilde{V}\rangle_H
   $$
and in particular
  \begin{equation} \label{dn}
  0 \le \|\widetilde{V}\|_H^2 =  \int_{\mathbb{S}^{d-1}} \langle \widetilde{V}(\rho \phi),\phi \rangle  \,\dd\sigma(\phi) =
   \int_{\mathbb{S}^{d-1}}  \int_{\mathbb{S}^{d-1}} \langle b(\rho \theta-\rho \phi)\theta,\phi \rangle\,\dd\sigma(\theta)\dd\sigma( \phi).
 \end{equation}
If the right side of (\ref{dn}) is non-zero then taking $V =
-\widetilde{V}$ gives
   $$
 \langle V(\rho \theta),\theta\rangle = -\langle \widetilde{V}(\rho \theta),\theta\rangle
 = -  \int_{\mathbb{S}^{d-1}} \langle \widetilde{V}(\rho \theta),\theta \rangle
 \dd\sigma(\theta) =- \|\widetilde{V}\|^2_H < 0,
   $$
so that (\ref{Vin}) is satisfied.  Conversely, if the right side
of (\ref{dn}) is zero then $\|\widetilde{V}\|_H = 0$ and so
   $$
 \int_{\mathbb{S}^{d-1}} \langle V(\rho \theta),\theta \rangle
 \dd\sigma(\theta) = \langle V, \widetilde{V}\rangle_H = 0
   $$
for all $V \in H$, so that $V$ cannot satisfy (\ref{Vin}).
Therefore, in order to prove Theorem \ref{thm-V}, it remains to
show that the right side of (\ref{dn}) is positive if and only if
$\mathbf{(C)_\rho}$ is satisfied.   The proof of Theorem
\ref{thm-V} is completed using the following result.

\begin{propo}  \label{prop-Bes} For $\rho > 0$
 \begin{equation} \label{Bes}
    \int_{\mathbb{S}^{d-1}} \int_{\mathbb{S}^{d-1}} \langle b(\rho \theta- \rho\phi)\theta,\phi\rangle \,\dd\sigma(\theta)\dd\sigma(\phi)
           = \mu_1 2^{d-2}\left[\Gamma\left(\frac{d}{2}\right)\right]^2 \int_{(0,\infty)} \left(\frac{{\cal J}_{d/2}(\rho s)
       }{(\rho s)^{(d-2)/2}}\right)^2 \dd M_P(s).
 \end{equation}
In particular
   $$
    \int_{\mathbb{S}^{d-1}} \int_{\mathbb{S}^{d-1}} \langle b(\rho \theta- \rho\phi)\theta,\phi\rangle \,\dd\sigma(\theta)\dd\sigma(\phi)
    > 0
   $$
if and only if condition $\mathbf{(C)_\rho}$ is satisfied.
    \end{propo}

\noindent{\bf Proof.}  We consider three special cases
corresponding to the three terms in the decomposition
(\ref{pre_ibf_decom}).  First, if $b_{jk}(x) = \delta_{jk}$, then
     $$
     \int_{\mathbb{S}^{d-1}} \int_{\mathbb{S}^{d-1}} \langle b(\rho \theta- \rho\phi)\theta,\phi\rangle \,\dd\sigma(\theta)\dd\sigma(\phi)
     = \int_{\mathbb{S}^{d-1}}\int_{\mathbb{S}^{d-1}} \langle \theta, \phi \rangle  \,\dd\sigma(\theta)\dd\sigma(\phi)
     = 0.
     $$
For the second case we suppose that $b_{jk}(x)$ is isotropic
potential.  In this case equation (4.25) of Yaglom \cite{Yag57}
gives
   $$
   b_{jk}(x) = \int_{(0,\infty)}  \int_{\mathbb{S}^{d-1}} e^{i\langle x,s \psi\rangle }\psi_j \psi_k \dd\sigma(\psi)\dd M_P(s)
    $$
and so
  $$
  \langle b(\rho \theta- \rho \phi)\theta,\phi\rangle
   = \int_0^\infty\int_{\mathbb{S}^{d-1}} \langle \theta,\psi\rangle
    \langle \phi, \psi \rangle e^{i\langle \rho \theta- \rho \phi, s \psi\rangle } \dd\sigma(\psi)
    \dd M_P(s).
   $$
Applying Fubini's theorem to the bounded function $\langle
\theta,\psi\rangle
    \langle \phi, \psi \rangle e^{i\langle \rho \theta- \rho \phi, s \psi\rangle
    }$
we get
   \begin{eqnarray*}
   \lefteqn{
  \int_{\mathbb{S}^{d-1}} \int_{\mathbb{S}^{d-1}} \langle b(\rho
  \theta-\rho
  \phi)\theta,\phi \rangle \,\dd\sigma(\theta)\dd\sigma(\phi)} \\
   &  = &  \int_{\mathbb{S}^{d-1}} \int_{\mathbb{S}^{d-1}}\left(\int_0^\infty \int_{\mathbb{S}^{d-1}} \langle \theta,\psi\rangle
    \langle \phi, \psi \rangle e^{i\langle \rho \theta-\rho \psi,s \psi\rangle }  \dd\sigma(\psi)
    \dd M_P(s) \right) \dd\sigma(\theta)\dd\sigma(\phi)\\
     &  = &  \int_0^\infty\int_{\mathbb{S}^{d-1}}\left( \int_{\mathbb{S}^{d-1}} \int_{\mathbb{S}^{d-1}} \langle \theta,\psi\rangle
    \langle \phi, \psi \rangle e^{i\langle \rho \theta- \rho \psi,s \psi\rangle }\dd\sigma(\theta)\dd\sigma(\phi)\right) \dd\sigma(\psi)
    \dd M_P(s)\\
         &  = &  \int_0^\infty\int_{\mathbb{S}^{d-1}} \left(\int_{\mathbb{S}^{d-1}}  \langle \theta,\psi\rangle
    e^{i\langle \rho \theta, s \psi\rangle } \dd\sigma(\theta)\right)\left(\int_{\mathbb{S}^{d-1}}\langle \phi, \psi \rangle
     e^{-i\langle \rho \psi,s \psi\rangle }\dd\sigma(\phi)\right) \dd\sigma(\psi)
    \dd M_P(s).
  \end{eqnarray*}
Now
  \begin{eqnarray*}
   \int_{\mathbb{S}^{d-1}}  \langle \theta,\psi\rangle
    e^{i\langle \rho \theta, s \psi\rangle } \dd\sigma(\theta)
    &  = &
     \frac{1}{is} \frac{\partial}{\partial \rho} \int_{\mathbb{S}^{d-1}}  e^{i\langle \rho \theta,s \psi\rangle }
     \dd\sigma(\theta)\\[1ex]
       & \stackrel{(i)}{=} & \frac{2^{(d-2)/2}\Gamma(\frac{d}{2})}{is} \frac{\partial}{\partial \rho} \left(\frac{{\cal J}_{(d-2)/2}(\rho s)
       }{(\rho s)^{(d-2)/2}}\right)\\[1ex]
       & \stackrel{(ii)}{=} & \frac{-2^{(d-2)/2}\Gamma(\frac{d}{2})}{i}  \left(\frac{{\cal J}_{d/2}(\rho s)
       }{(\rho s)^{(d-2)/2}}\right)
  \end{eqnarray*}
and similarly
   \begin{eqnarray*}
   \int_{\mathbb{S}^{d-1}}  \langle \phi,\psi\rangle
    e^{-i\langle \rho \phi,s \psi\rangle } \dd\sigma(\phi)
    &  = &
     \frac{2^{(d-2)/2}\Gamma(\frac{d}{2})}{i}  \left(\frac{{\cal J}_{d/2}(\rho s)
       }{(\rho s)^{(d-2)/2}}\right).
  \end{eqnarray*}
The equality ({\em i}) uses the fact that
 $$
  \int_{\mathbb{S}^{d-1}} e^{i \langle x,  \theta\rangle}  \dd\sigma(\theta)
   = \frac{\Gamma(\frac{d}{2}) {\cal J}_{(d-2)/2}(|x|) }{(|x|/2)^{(d-2)/2}},
    $$
see equations (4.9) through (4.11) of Yaglom \cite{Yag57}.  For
({\em ii}) see, for example, Watson \cite[page 45]{Wat44}.
Therefore
   \begin{eqnarray*}
   \lefteqn{
  \int_{\mathbb{S}^{d-1}} \int_{\mathbb{S}^{d-1}} \langle b(\rho \theta-\rho \phi)\theta,\phi\rangle \,\dd\sigma(\theta)\dd\sigma(\phi)} \\
      &  = &  \int_0^\infty\int_{\mathbb{S}^{d-1}} \left(\int_{\mathbb{S}^{d-1}}  \langle \theta,\psi\rangle
    e^{i\langle \rho \theta,s \psi\rangle } \dd\sigma(\theta)\right)
    \left(\int_{\mathbb{S}^{d-1}}\langle \phi, \psi \rangle e^{-i\langle \rho \psi, s \psi\rangle }\dd\sigma(\phi)\right) \dd\sigma(\psi)
    \dd M_P(s)\\
     & = &  2^{d-2}\left[\Gamma(d/2)\right]^2 \int_0^\infty
    \left(\frac{{\cal J}_{d/2}(\rho s)
       }{(\rho s)^{(d-2)/2}}\right)^2 \dd M_P(s).
      \end{eqnarray*}
For the third case we suppose that $b_{jk}(x)$ is isotropic
solenoidal, so that $\mbox{div}V = 0$ for every $V \in H$. For the
vector field $\widetilde{V}$ defined in (\ref{D}) the divergence
theorem applied to the ball of radius $\rho$ gives
$$
   \langle V,\widetilde{V}\rangle = \int_{\mathbb{S}^{d-1}} \langle V(\rho \theta), \theta \rangle
   \,dm(\theta) = 0
   $$
for all $V \in H$ because $\mbox{div}\widetilde{V} = 0$.  Taking
$V = \widetilde{V}$ we get
   $$
  \int_{\mathbb{S}^{d-1}}\int_{\mathbb{S}^{d-1}} \langle b(\rho \theta-\rho \phi)\theta,\phi \rangle \,d\sigma(\theta)d\sigma(\phi)
  = \|\widetilde{V}\|^2 = 0.
  $$
The first assertion is now obtained by summing the contributions
of the three terms in the decomposition (\ref{pre_ibf_decom}). The
second assertion follows because the right side of (\ref{Bes}) is
strictly positive unless ${\cal J}_{d/2}(\rho s) = 0$ for
$M_P$-almost all $s  \in (0,\infty)$. \hfill $\Box$

\section{A control theorem for stochastic flows}
\label{sec-control}

The following theorem is based on a result of Dolgopyat, Kaloshin
and Koralov \cite[Section 2.4]{DKK04a}.  It does not use the
isotropy condition.

\begin{theo}\label{sb_thm} Assume that $ v(t,x): \R^d\times \R_+
\to \R^d $ is a deterministic function which is $C^2$ with bounded
derivatives in the space variable $x$ for every $t$ and continuous
in $t$ for every $x$.  Let $M(t,x)$ be a generating Brownian field
with covariance tensor $b(x,y)$ which is $C^{2,2}$ with bounded
derivatives in the variables $x$ and $y$.  Let $\{\phi_{s,t}(x): 0
\le s \le t < \infty,\, x \in \IR^d \}$ be the flow generated by
the semimartingale field $F(t,x)=M(t,x)+\int_0^t v(s,x)\dd s$,
i.e.\
$$\phi_{s,t}(x)=x+\int\limits_s^tF(\dd u,\phi_{s,u}(x))\,\,\, \text{ for all }x\in\R^d, \,\,0\le s \le t\,. $$
Suppose that $V$ is a vector field in the reproducing kernel
Hilbert space of $M$, and let $\xi_t(x)$ denote the solution of
the ordinary differential equation
 $$
 \left\{ \begin{array}{rcl}
\dot{\xi}_t(x) & = & V(\xi_t(x)), \\
\xi_0(x) & =& x.
\end{array} \right .
  $$
Then for all compact subsets $K\subset \R^d$, time instants $T>0$
and all $\delta>0$ there exist $c_0 > 0$ such that
\begin{align*}
\P\big(\,\,\sup_{t \le T} \sup _{x \in K}\big|
\phi_{t/c}(x)-\xi_t(x)\big| <\delta \big)>0
\end{align*}
for all $c \ge c_0$.
\end{theo}

\noindent \textbf{Proof}: Choose a compact set $K\subset \R^d$, a
time $T>0$ and $\delta>0$.  As usual we will write $ \phi_t $ for
$\phi_{0,t}$. Let $\{V_i: i \ge 1\}$ be a complete orthonormal set
of vector fields for the reproducing kernel Hilbert space $H$.
Then as in \cite{BaH86} (see also Bogachev \cite[Theorem
3.5.1]{Bog98}), the generating Brownian field $M(t,x)$ can be
written in distribution as
  $$
  M(t,x) = \sum_{i = 1}^\infty B^i(t)V_i(x)
  $$
where $\{B^i(t): t \ge 0 \}_{i \ge 1}$ are independent standard
scalar Brownian motions, and the Kunita stochastic differential
equation for $\phi$ can be written
  \begin{align}\label{sb_ct_sde}
\phi_t(x)=x+\sum_{i=1}^\infty\int\limits_0^t V_i(\phi_s(x))\dd
B_s^i+\int\limits_0^tv(s,\phi_s(x))\dd s\,.
  \end{align}
We will work with a complete orthonormal set chosen so that $V =
\lambda V_1$ for some $\lambda > 0$.  For any $c > 0$ we will
consider the SDE (\ref{sb_ct_sde}) for $0 \le t \le T/c$.  Without
loss of generality we can assume that the underlying probability
space $(\Omega,\F,\P)$ is the infinite product of the path spaces
of the standard Brownian motions $B^i$, that is
\begin{align*}
(\Omega,\F,\P) = \left(\underset{i\in \N}  \times C([0,T/c]:\R)
,\underset{i\in \N}\otimes \B(C([0,T/c]:\R)),\underset{i\in
\N}\otimes \P_i\right),
\end{align*}
where $\P_i$ denotes the Wiener measure on
$(C([0,T/c]:\R),\B(C([0,T/c]:\R)))$ for all $i \in \N$.  On the
canonical path space we have $B_t^i(\omega) = \omega^i(t)$ for all
$\omega \in \Omega$ and $t\le T/c$.
  Define $B_t^c=
(B_t^{c,1},B_t^{c,2},\dots):=(B_t^1- \lambda c t, B_t^2, \dots)$.
Then the SDE (\ref{sb_ct_sde}) becomes
\begin{align*}
 \phi_t(x)=x+c\int\limits_0^t V (\phi_s(x)) \dd s+\sum_{i=1}^{\infty}\int\limits_0^t V_i(\phi_s(x))\dd B_s^{c,i}
    +\int\limits_0^tv(s,\phi_s(x))\dd s
    \end{align*}
for $t \in [0,T/c]$.  Let $\P_1^c$ be the measure on
$(C([0,T/c]:\R),\B(C([0,T/c]:\R)))$ under which $B^{c,1}$ is a
standard Brownian motion.  It is well known that $\P_1^c$ is
equivalent to $\P_1$ with
   $$
  \P_1^c(\dd\omega^1)=\exp\left[\lambda c \,\omega^1(T/c)-\frac{\lambda^2
c T}{2}\right] \P_1(\dd\omega^1)\,.
    $$
  Since
the transformation of $B$ into $B^c$ involves only the first
coordinate, it follows that $B^c$ is a standard
infinite-dimensional Brownian motion under $\P^c:=\P_1^c \otimes
\P_2\otimes\P_3\otimes \dots$ and that $\P^c$ is equivalent to
$\P$.

Consider the time changed process $t \mapsto \phi_{t/c}(x)$. For
$0 \le t \le T$ we have
\begin{eqnarray*}
 \phi_{t/c}(x)
  & = &  x+c\int\limits_0^{t/c} V(\phi_s(x)) \dd s+\sum_{i=1}^{\infty}\int\limits_0^{t/c} V_i(\phi_s(x))\dd B_s^{c,i}
   + \int\limits_0^{t/c} v(s, \phi_s(x))\dd s\\
 & = & x +\int\limits_0^t V(\phi_\frac{s}{c}(x)) \dd s + \frac{1}{\sqrt{c}}\sum_{i=1}^{\infty}\int\limits_0^t
   V_i(\phi_{s/c}(x))\dd \widetilde{B}_s^i
      +\frac{1}{c}\int\limits_0^t v(s/c,\phi_{s/c}(x))\dd s
\end{eqnarray*}
where $\widetilde{B}_t^i = \sqrt{c} B^{c,i}_{t/c}$ for $0 \le t
\le T$ and $i \ge 1$.  (For the time change in the stochastic
integral see for example {\O}ksendal \cite[Theorem 8.20]{Oks03}.)
Under $\P^c$, the process $\widetilde{B} = \{\widetilde{B}_t^i: 0
\le t \le T\}_{i \ge 1}$ is a standard infinite-dimensional
Brownian motion.  It follows that the distribution of
$\left\{\phi_{t/c}(x) \colon t \in [0, T], x\in \R^d \right\}$
under the probability measure $\P^c$ is the same as the
distribution of $\left\{Y^c_t(x)\colon t \in [0,T], x\in \R^d
\right\}$ generated by the SDE
\begin{align}\label{sb_ct2}
Y_t^c(x)=x+\int\limits_0^t V(Y_s^c(x)) \dd s+\frac{1}{\sqrt
c}\sum_{i=1}^{\infty}\int\limits_0^t V_i(Y_s^c(x))\dd\tilde{B}_s^i
+\frac{1}{c}\int\limits_0^tv({s/c},Y_s^c(x))\dd s\,,
\end{align}
for $ t \in [0,T]$ and $ x\in \R^d $, where $ \widetilde{B} =
\{\widetilde{B}_t^i: 0 \le t \le T\}_{i \ge 1}$ is a standard
infinite dimensional Brownian motion on some probability space
$(\widetilde \Omega,\widetilde \F, \widetilde \P)$. Since
  \begin{equation} \label{diff}
  Y_t^c(x)-\xi_t(x) =
  \int_0^t\left(V(Y_s^c(x))-V(\xi_s(x))\right)\dd s  +\frac{1}{\sqrt
c}\sum_{i=1}^{\infty}\int_0^t V_i(Y_s^c(x))\dd \widetilde{B}_s^i
+\frac{1}{c}\int_0^t v({s/c},Y_s^c(x))\dd s
  \end{equation}
for any $p \ge 1$ there exists $L_p$ such that
   $$
   \widetilde{\E}\sup_{t \le T} |Y_t^c(x)-\xi_t(x)|^p \le L_p\left(\frac{1}{c^p} + \frac{1}{c^{p/2}}\right)
   $$
for all $x \in \mathbb{R}^d$.  Similarly, taking derivatives with
respect to $x$ in (\ref{diff}), there exists $M_p$ such that
   $$
  \widetilde{\E} \sup_{t \le T} \left\|DY_t^c(x) - D\xi_t(x)\right\|^p \le M_p\left(\frac{1}{c^p} + \frac{1}{c^{p/2}}\right)
  $$
for all $x \in \R^d$.  Taking $p > d$ and using the Sobolev
embedding theorem, for any compact set $K$ there exists a constant
$N_p$ such that
      \begin{equation} \label{sb_ct4}
   \widetilde{\E}\left(\sup_{t \le T} \sup_{x \in K}|Y_t^c(x)-\xi_t(x)|^p\right) \le N_p \left(\frac{1}{c^p} +
   \frac{1}{c^{p/2}}\right).
    \end{equation}
Details of similar calculations can be found in Kunita
\cite[Section 5.4]{Ku90} and Ikeda and Watanabe \cite[Section
V.2]{IW81}.  Therefore there is $c_0 > 0$ such that if $c \ge c_0$
then
\begin{align*}
\widetilde \P\big(\,\,\sup_{t \le T}\sup_{x \in
K}\big|Y^c_t(x)-\xi_t(x)\big|< \delta\big) \ge \frac{1}{2}\,,
\end{align*}
and since the distribution of $\big\{Y^c_t(x):x \in \R^d, t\in
[0,T]\big\}$ under $\widetilde\P$ coincides with the distribution
of $\big\{\phi_{t/c}(x):x \in \R^d, t\in [0,T]\big\} $ under
$\P^c$, we obtain
\begin{align*}
\P^c\big(\,\,\sup_{t \le T}\sup_{x \in
K}\big|\phi_{{t/c}}(x)-\xi_t(x)\big|< \delta\big) \ge
\frac{1}{2}\,.
\end{align*}
Since $\P^c \sim\P$ we obtain
\begin{align*}
\P\big(\,\,\sup_{t \le T}\sup_{x \in
K}\big|\phi_{{t/c}}(x)-\xi_t(x)\big|< \delta\big) > 0\,,
\end{align*}
which is the statement of the theorem. \hfill $\square$

\section{Proofs for Section \ref{sec-main}} \label{sec-pf}

\textbf{Proof of Theorem \ref{mt_thm1}}:  According to Theorem
\ref{thm-V}, under condition $\mathbf{(C)}_R$ the reproducing
kernel Hilbert space $H$ contains a rotation invariant vector
field $V$ pointing inwards everywhere on the surface of the ball
$\text{B}(0,R)$, that
 is,
there is a constant $\alpha>0$ such that
\begin{align*}
\langle V(x),x \rangle \le -\alpha<0 \quad \textrm{ for all } x
\in\partial \text{B}(0,R).
\end{align*}
The continuity of $V$ implies the existence of a $\delta>0$ such
that
\begin{align}\label{sb_thm2_1}
\langle V(x),x \rangle \le -\frac{\alpha}{2}<0 \quad \textrm{ for
all } x \in  \text{B}(0,R+\delta)\setminus
\text{B}(0,R-3\delta)\,.
\end{align}
Let $\xi_t$ denote the time $t$ flow along the vector field $V$,
as in Theorem \ref{sb_thm}.  Using the condition
(\ref{sb_thm2_1}), there is a time $T_3>0$ such that
\begin{equation}\label{mt_thm1_eq1}
 \left\{
   \begin{array}{ll}
   \xi_t(\text{B}(0, R-2\delta)) \subset \text{B}(0, R-2\delta) & \mbox{ for
   } 0 \le t \le T_3, \\[1ex]
   \xi_t(\text{B}(0, R+\delta)) \subset \text{B}(0, R-3\delta) & \mbox{ for } t \ge T_3.
   \end{array} \right.
\end{equation}
  Now apply the Theorem  \ref{sb_thm} with $T =
2T_3$ and some compact set $K$ containing $\text{B}(0,R+\delta)$,
so that there exists $c \ge T_3/T_1$ for which
 \begin{align}\label{mt_thm1_eq2}
\P\left(\,\sup_{0 \le t \le 2T_3}\sup_{x \in
\mathrm{B}(0,R+\delta)}\big|\phi_{t/c}(x)-\xi_t(x)\big|< \delta
\right)> 0 \,.
\end{align}
Using (\ref{mt_thm1_eq1}) we have
   \begin{eqnarray*}
   \lefteqn{ \left\{\,\sup_{0 \le t \le 2T_3}\sup_{x \in
\mathrm{B}(0,R+\delta)}\big|\phi_{t/c}(x)-\xi_t(x)\big|< \delta \right\}} \hspace{4em}
\\[1ex]
          & \subset &
     \left\{\,    \phi_{t/c}(\mathrm{B}(0,R-2\delta)) \subset \mathrm{B}(0,R-\delta) \mbox{ \rm{for} } 0 \le t \le T_3
     \right\}\\[0.5ex]
     & & \hspace{1em}
      \cap  \left\{\,    \phi_{t/c}(\mathrm{B}(0,R+\delta)) \subset \mathrm{B}(0,R-2\delta) \mbox{ \rm{for} } T_3 \le t \le 2T_3 \right\}
      \\[1ex]
     & = &
     \left\{\,    \phi_{t}(\mathrm{B}(0,R-2\delta)) \subset \mathrm{B}(0,R-\delta) \mbox{ \rm{for} } 0 \le t \le T_3/c
     \right\}\\[0.5ex]
     & & \hspace{1em}
      \cap
      \left\{\,    \phi_{t}(\mathrm{B}(0,R+\delta)) \subset \mathrm{B}(0,R-2\delta) \mbox{ \rm{for} } T_3/c \le t \le 2T_3/c
      \right\}\\[1ex]
     & \equiv & A_0,
     \end{eqnarray*}
say, and then (\ref{mt_thm1_eq2}) gives:
   \begin{equation} \label{mt_thm1_eq3}
     \P(A_0) \ge  \P\left(\,\sup_{t \le 2T_3}\sup_{x \in
K}\big|\phi_{t/c}(x)-\xi_t(x)\big|< \delta \right)> 0.
   \end{equation}
Let $A_k$ denote the time $2kT_3/c$ shift of the event $A_0$, so
that
  \begin{eqnarray*}
     A_k & = &
       \left\{\,    \phi_{2kT_3/c,t}(\mathrm{B}(0,R-2\delta)) \subset \mathrm{B}(0,R-\delta) \mbox{ \rm{for} } 2kT_3/c \le t \le (2k+1)T_3/c
     \right\}\\[1ex]
     & & \hspace{1em}
      \cap
      \left\{\,    \phi_{2kT_3/c,t}(\mathrm{B}(0,R+\delta)) \subset \mathrm{B}(0,R-2\delta) \mbox{ \rm{for} } (2k+1)T_3/c \le t \le (2k+2)T_3/c
      \right\}.
       \end{eqnarray*}
The events $\{A_k: k \ge 0\}$ are independent with $P(A_k) =
P(A_0)$ for all $k \ge 1$, and
   $$
   \bigcap_{k=0}^{n-1} A_k \subset \left\{\,    \phi_{t}(\mathrm{B}(0,R+\delta)) \subset \mathrm{B}(0,R-\delta) \mbox{ \rm{for} } T_3/c \le t \le 2nT_3/c
     \right\}.
     $$
Therefore
 \begin{eqnarray*}
  \lefteqn{\P\left(  \phi_{t}(\mathrm{B}(0,R+\delta)) \subset \mathrm{B}(0,R-\delta) \mbox{ \rm{for} } T_3/c \le t \le 2nT_3/c
     \right)} \hspace{12em}\\[2ex]
     &  \ge &  \P\left(\bigcap_{k=0}^{n-1} A_k\right) =
     \prod_{k=0}^{n-1}\P(A_k) = \left(\P(A_0)\right)^n > 0.
     \end{eqnarray*}
Since $T_3/c \le T_1$ and $n$ can be chosen so that $2nT_3/c \ge
T_2$, the proof of the squeezing result (\ref{Rcont}) is complete.
The corresponding expansion result (\ref{Rexp}) can be proved
using the same method.  Simply replace $\xi_t$ by the time $t$
flow $\widehat{\xi}_t$ along $-V$, so that (\ref{mt_thm1_eq1}) can
be replaced by
  $$
 \left\{
   \begin{array}{ll}
   \widehat{\xi}_t(\text{B}(0, R+2\delta)) \supset \text{B}(0, R+2\delta) & \mbox{ for
   } 0 \le t \le T_3, \\[1ex]
   \widehat{\xi}_t(\text{B}(0, R-\delta)) \supset \text{B}(0, R+3\delta) & \mbox{ for } t \ge
   T_3,
   \end{array} \right.
  $$
and then $\P(\widehat{A}_0) > 0$ where
    \begin{eqnarray*}
     \widehat{A}_0 & = &
       \left\{\,    \phi_{t}(\mathrm{B}(0,R+2\delta)) \supset \mathrm{B}(0,R+\delta) \mbox{ \rm{for} } 0 \le t \le T_3/c
     \right\}\\[1ex]
     & & \hspace{1em}
      \cap
      \left\{\,    \phi_{t}(\mathrm{B}(0,R-\delta)) \supset \mathrm{B}(0,R+2\delta) \mbox{ \rm{for} } T_3/c \le t \le 2T_3/c
      \right\}.
       \end{eqnarray*}
\hfill$\square$

\vspace{2ex}

 \noindent \textbf{Proof of  Theorem \ref{mt_thm2}}:
According to (\ref{Rcont}) of Theorem \ref{mt_thm1}, for each $u
\in [r,R]$ there exists $\delta(u)>0$ such that
\begin{align*}
\P\big[\phi_{t}\left(\mathrm{B}(0,u+\delta(u))\right)\subset
\mathrm{B}(0,u-\delta(u)) \mbox{ for  } t_1 \le t \le t_2  \big]
>0\,
\end{align*}
for any $0 < t_1 < t_2$.  Consider the open cover
$$\{(u-\delta(u),u+\delta(u)) \,\colon\, u \in [r,R] \}$$ of the
compact interval $[r,R]$, and choose a minimal finite subcover
 $$\{(u_i-\delta(u_i),u_i+\delta(u_i)) \,\colon\, i=1,\dots,N \}\,.$$
Under the condition of minimality, the intervals can be labelled
so that $u_1+ \delta(u_1) > R$, $u_N-\delta(u_N) < r$ and
$u_i-\delta(u_i)<u_{i+1}+\delta(u_{i+1})<u_i+\delta(u_i)$ for all
$i \in \{1, \dots , N-1\}$.
 Set $\overline  u_i:=u_{i}+ \delta(u_i)$ and  $\underline u_i=u_{i}-
\delta(u_i)$. With this notation, the conditions on the subcover
become:  $\overline u_1 \ge R$, $\underline u_N \le r$ and
$\underline u_i< \overline  u_{i+1}<\overline  u_i$ for all $i \in
\{1, \dots , N-1\}$. Choose times $0 < t_1< \cdots < t_{N-1} <
T_1$ and consider the events
   $$
   A_i = \left\{\phi_{t_{i-1},t_i}(\mathrm{B}(0,\overline u_i)) \subset
   \mathrm{B}(0,\underline u_i) \right\}
   $$
for $1 \le i \le N-1$, and
  $$
   A_N =  \left\{\phi_{t_{N-1},t}(\mathrm{B}(0,\overline u_N)) \subset
   \mathrm{B}(0,\underline u_N)\mbox{ for } T_1 \le t \le T_2 \right\}
   $$
The events $\{A_k: 1 \le k \le N\}$ are independent, and by
Theorem \ref{mt_thm1} we have $\P(A_k) > 0$ for $1 \le k \le N$.
Moreover
  \begin{eqnarray*}
\bigcap_{i=1}^N A_k
        & \subset &
       \bigcap_{i=1}^{N-1} \big\{
        \phi_{t_{i-1},t_i}\left(\mathrm{B}(0,\overline u_i)\right) \subset \mathrm{B}(0,\overline u_{i+1})
        \big\} \\
        & & \hspace{4ex} \cap \,\big\{
        \phi_{t_{N-1},t}\left(\mathrm{B}(0,\overline u_N)\right) \subset \mathrm{B}(0,\underline
        u_N)\mbox{ for } T_1 \le t\le T_2
        \big\}\\[2ex]
        & \subset &
           \big\{
        \phi_{0,t}\left(\mathrm{B}(0,\overline u_1)\right) \subset \mathrm{B}(0,\underline
        u_N)\mbox{ for } T_1 \le t\le T_2
        \big\}\\[2ex]
        & \subset &
           \big\{
        \phi_{0,t}\left(\overline{\mathrm{B}}(0, R)\right) \subset
        \mathrm{B}(0,r)\mbox{ for } T_1 \le t \le T_2
        \big\}.
         \end{eqnarray*}
Therefore
  $$
    \P\big(
        \phi_{0,t}\left(\overline{\mathrm{B}}(0, R)\right) \subset
        \mathrm{B}(0,r)\mbox{ for } T_1 \le t \le T_2
        \big) \ge \P\left(\bigcap_{k=1}^N A_k \right) = \prod_{k=1}^N \P(A_k) > 0
        $$
and the proof of the squeezing result is complete.  The proof of
the expanding result uses the same method with (\ref{Rexp}) in
place of (\ref{Rcont}).  \hfill $\square$

\vspace{2ex}

\begin{remark}{\rm  According to Theorem \ref{thm-V}, if
the condition $\mathbf{(C)_\rho}$ is not satisfied there is no
vector field $V$ in the RKHS pointing strictly inwards on the
boundary $\partial B(0,\rho)$ and therefore we cannot apply our
technique to show the ball of radius $\rho $ can be squeezed.  If
$\mu_1 > 0$ it is an open question whether there exist two (or
more) vector fields $V_1$, $V_2$ from the RKHS such that when
applying them consecutively to $B(0,\rho)$ the resulting image is
strictly contained in $B(0,\rho)$.  However if $\mu_1 = 0$ then
the isotropic Brownian field is almost surely divergence-free, so
that the resulting IBF is volume preserving and the conclusion of
Theorem \ref{mt_thm1} is definitely false.}
 \end{remark}

\section{Lengths of Curves} \label{sec-length}

For fixed $x \in \IR^d$ and $v \in \IR^d$ the asymptotic rate of
growth or decay of $v_t = D\phi_t(x)v$ is given by the top
Lyapunov exponent
   $$
   \lambda = \lim_{t \to \infty} \frac{1}{t} \log |v_t| \quad
   \quad \P \mbox{ almost surely}.
   $$
For an IBF the limit exists and takes the same value for all $x$
and $v$.  It is shown in \cite{BaH86} and \cite{LeJ85} that
   $$
   \lambda = (d-1)\frac{\beta_N}{2}-\frac{\beta_L}{2} $$
where $\beta_L$ and $\beta_N$ are given in (\ref{betaL}) and
(\ref{betaN}).

The top Lyapunov exponent describes the rate of growth or decay of
a single tangent tangent vector.  Here we will consider the more
complicated situation concerning the growth or decay of the length
of a differentiable curve under an IBF.  Throughout this section
we will assume that $\Gamma$ is a piecewise $C^1$ curve in $\R^d$
parameterized by $\gamma(u):[0,1]\to \R^d$ with $\gamma'(u) \ne 0$
at all $C^1$ points $u \in [0,1]$. The image $\phi_t(\Gamma)$ of
$\Gamma$ under the action of the flow will be denoted by
$\Gamma_t$. It is a piecewise $C^1 $ curve parameterized by
$\gamma_t = \phi_t \circ \gamma:[0,1] \to \R^d$.  The length of
$\Gamma_t$ will be denoted by $L_t$ and is given by
  $$
  L_t = \int\limits_0^1 |
\gamma_t'(u)|\dd u = \int_0^1 |D\phi_t(\gamma(u))\gamma'(u)|\dd u.
 $$
The diameter of the curve $\Gamma_t$ is given by
$\mbox{diam}(\Gamma_t) = \sup\{|x-y|: x,y \in \Gamma_t\}$.

\vspace{2ex}

It has been shown, see Dimitroff \cite[Proposition 3.1.2]{Dim06},
that the evolving length of a differentiable curve in an IBF
satisfies the following bounds

\begin{align*}
\lambda \le \liminf_{t\to \infty} \frac{1}{t}\ln L_t \le
\limsup_{t\to \infty} \frac{1}{t}\ln L_t \le
\lambda+\frac{\beta_L}{2}.
\end{align*}
It follows that if the top Lyapunov exponent $\lambda$ is positive
then the length exhibits exponential growth.  However if $\lambda
< 0$ the above bounds are not very informative because
$\lambda+\beta_L/2 = (d-1)\beta_N/2$ is always strictly positive.
Here we will be dealing with the case $\lambda < 0$.  Baxendale
and Harris have shown in \cite{BaH86} that if $d=2$ and
$\frac{\beta_L}{\beta_N}>\frac{5}{3}$ then the length of a curve
decreases to $0$ on the set where the diameter of the curve
decreases to $0$, provided of course this happens with positive
probability, i.e.
\begin{align}\label{le_mstat}
\P( L_t \to 0 \mbox{ as } t \to \infty \,|\,
\text{diam}(\Gamma_t)\to 0 \mbox{ as } t \to \infty)=1
\end{align}
provided that $\P( \text{diam}(\Gamma_t)\to 0 \mbox{ as } t \to
\infty)>0$.  The following result of Dimitroff \cite[Theorem
3.2.1]{Dim06} strengthens (\ref{le_mstat}) by giving an
exponential decay rate under more general conditions.  (In
dimension $d=2$ the condition $\lambda < 0$ is equivalent to
$\frac{\beta_L}{\beta_N}>1$.)

\begin{theo}\label{le_shrink_main}
Let $\{\phi_{s,t}(x): 0 \le s \le t < \infty,\,x \in \IR^d \}$ be
a $d$-dimensional isotropic Brownian flow with top Lyapunov
exponent $\lambda < 0$. Then
\begin{align*}
\P\left( \left.\lim_{t \to \infty}\frac{1}{t}\log L_t=\lambda
\,\,\right|  \,\, \textnormal{diam}(\Gamma_t) \to 0 \mbox{ as }t
\to \infty  \right)=1
\end{align*}
provided $\P\left(\textnormal{diam}(\Gamma_t) \to 0 \mbox{ as } t
\to \infty\right)>0$.
\end{theo}

We can apply our Corollary \ref{cor1} to show that
$\P\left(\textnormal{diam}(\Gamma_t) \to 0 \mbox{ as } t \to
\infty\right)>0$ provided that condition $\mathbf{(C)_\rho}$ holds
for all $\rho > 0$.  This will be done in Proposition
\ref{ap_curve} below.  Before doing so we first sketch the proof
of the above theorem.  Concepts and notation from this proof will
then be used in the proof of Proposition \ref{ap_curve}.

\vspace{2ex}

\noindent \textbf{Sketch of the proof:}  The detailed proof of
Theorem \ref{le_shrink_main} can be found in \cite[Theorem
3.2.1]{Dim06}.  Since the distribution of the IBF is translation
invariant we can assume without loss of generality that $0 \in
\Gamma$.  We first consider the flow
$\psi_{s,t}(x):=\phi_{s,t}(x+\phi_{0,s}(0))-\phi_{0,t}(0)\,,$
which is simply the original flow viewed from the point of view of
the moving particle started in the origin.  Then $\psi$ and $\phi$
have the same Lyapunov spectrum and the length of $\phi_t(\Gamma)=
\Gamma_t$ equals the length of $\psi_t(\Gamma)\equiv
\widehat{\Gamma}_t$. The flow $\psi$ can be extended to a double
sided random dynamical system with fixed point at the origin over
a measurable ergodic shift $\{\theta_t: t \in \IR\}$, and then we
can write $\psi_{s,t}(x,\omega) = \psi(t-s,x,\theta_s \omega)$,
see \cite[Lemma 3.2.1]{Dim06}. This allows us to use the local
stable manifold theorem, see Mohammed and Scheutzow
\cite{MoSch99}.
 One of the characterizations of the local stable
manifold is
\begin{equation} \label{msrho}
\mathcal{S}(\omega)=\left\{x \in \overline{ \text B}
(0,\rho(\omega)) \,\,\big|\,\,|\psi(n,x,\omega)|\le
\beta(\omega)\exp(\lambda+\epsilon)n \textrm { for all } n \in \N
\right\},
\end{equation}
where $1>\beta(\omega)>\rho(\omega)>0$ and $\epsilon$ is some
fixed positive number between $0$ and $-\lambda$.
  Since $\lambda
< 0$ we have $\mathcal{S}(\omega) = \overline{\text
B}(0,\rho(\omega))$. Since $\{\rho \circ \theta_n\}_{n \in \N}$ is
a stationary sequence of positive functions, for each $\omega$
such that $\mbox{diam}(\widehat{\Gamma}_t(\omega))\to 0$ as $t \to
\infty$ there exists $n = n(\omega) \in \N$ such that
$\mbox{diam}(\widehat{\Gamma}_n(\omega)) \le \rho(\theta_n
\omega)$ and so $\widehat{\Gamma}_n(\omega) \subset
\mathcal{S}(\theta_n\omega)$.  Since
   $$
   L_{t+n}(\omega) \le L_n(\omega) \sup  \left
   \{
   \frac{|\psi(t,x_1,\theta_n\omega)-\psi(t,x_2,\theta_n\omega)|}{|x_1-x_2|} :
 x_1 \ne x_2,\,\,\, x_1,\,x_2 \in \widehat{\Gamma}_n(\omega) \right \}
 $$
the result now follows from the estimate
 \begin{align}\label{msb} \limsup_{t \to \infty}\frac{1}{t} \log
\left [ \sup  \left
   \{
   \frac{|\psi(t,x_1,\omega)-\psi(t,x_2,\omega)|}{|x_1-x_2|} :
 x_1 \ne x_2,\,\,\, x_1,\,x_2 \in \mathcal{S}(\omega) \right \} \right
]\le \lambda \,,
\end{align}
see  \cite[Theorem 3.1(ii)(b)]{MoSch99}, with $\omega$ replaced by
$\theta_n\omega$.
  \hfill $\Box$

\vspace{2ex}

The next proposition concerns the probability of the event
$\{\text{diam}(\Gamma_t)\to 0\}$. It is a direct consequence of
our  main result (Theorem \ref{mt_thm2}).
\begin{propo}\label{ap_curve}
Let $\{\phi_{s,t}(x): 0 \le s \le t < \infty,\,x \in \IR^d \}$ be
a $d$-dimensional isotropic Brownian flow with top Lyapunov
exponent $\lambda < 0$.   Assume that condition
$\mathbf{(C)}_\rho$ is satisfied for all $\rho > 0$. Then for all
initial curves $\Gamma$
$$\P\left(\textnormal{diam}(\Gamma_t) \to 0 \mbox{ as } t \to \infty \right) >0.$$
\end{propo}

\noindent{\bf Proof:}
 Let $r,R>0$ be such that $\P(\rho(\omega)>r)>0$, where $\rho(\omega)$ is the radius of the local stable manifold $\mathcal{S}(\omega)$
as in (\ref{msrho}), and $\Gamma \subset \overline{\text B}(0,R)$.
According to Theorem \ref{mt_thm2}, there exists $T>0$, such that
$\P (\phi_T(\overline{\text{B}}(0,R),\omega) \subset
\text{B}(0,r))>0$. Rewriting (\ref{msb}) in terms of the IBF
$\phi$ we have
  $$
   \limsup_{t \to \infty}\frac{1}{t} \log
\left [ \sup  \left
   \{
   \frac{|\phi(t,x_1,\omega)-\phi(t,x_2,\omega)|}{|x_1-x_2|} :
 x_1 \ne x_2,\,\,\, x_1,\,x_2 \in \mathcal{S}(\omega) \right \} \right
]\le \lambda
  $$
so that $\mbox{diam}(\phi_t(S(\omega),\omega)) \to 0$ as $t \to
\infty$.  Using this result with $\omega$ replaced by $\theta_t
\omega$ we obtain
  \begin{eqnarray*}
   \P\left(\mbox{diam}(\Gamma_t(\omega)) \to 0 \mbox{ as } t \to
   \infty\right)
   & \ge & \P\left(\Gamma_T(\omega) \in S(\theta_T \omega)\right)\\
    & \ge & \P\left(\Gamma_T(\omega) \subset {\text B}(0,r) \mbox{ and } {\text B}(0,r)\subset  S(\theta_T \omega)\right)\\
    & \ge & \P\left(\phi_T(\overline{\text B}(0,R),\omega) \subset {\text B}(0,r) \mbox{ and } \rho(\theta_T \omega) > r\right)
    \end{eqnarray*}
Observe that $\rho(\theta_T\omega)$ is measurable with respect to
the $\sigma$-algebra generated by the increments
$\{\phi_{T,t}(\cdot,\omega): t \ge T\}$ and hence it is
independent of $\phi_T( \cdot,\omega)$.  Therefore
\begin{eqnarray*}
   \P\left(\mbox{diam}(\Gamma_t(\omega)) \to 0 \mbox{ as } t \to
   \infty\right)
    & \ge & \P\left(\phi_T(\overline{ \text B}(0,R) \subset B(0,r)\right)\P \left( \rho(\theta_T \omega) > r \right)\\
    & \ge & \P\left(\phi_T(\overline{ \text B}(0,R) \subset B(0,r)\right)\P \left( \rho(\omega) > r
    \right)\\
    & > & 0,
    \end{eqnarray*}
which completes the proof.\hfill $ \square $

\bibliographystyle{abbrv}

\bibliography{biblio}

\begin{thebibliography}{10}

\bibitem{Bax76}
P.~Baxendale.
\newblock Gaussian measures on function spaces.
\newblock {\em American Journal of Mathematics}, 98(4):891--952, 1976.

\bibitem{BaH86}
P.~Baxendale and T.~E. Harris.
\newblock Isotropic stochastic flows.
\newblock {\em Ann. Probab.}, 14(2):1155--1179, 1986.

\bibitem{Bog98}
V.~I. Bogachev.
\newblock {\em Gaussian measures}, volume~62 of {\em Mathematical Surveys and
  Monographs}.
\newblock American Mathematical Society, Providence, RI, 1998.

\bibitem{CSS99}
M.~Cranston, M.~Scheutzow, and D.~Steinsaltz.
\newblock Linear expansion of isotropic {B}rownian flows.
\newblock {\em Elect.~Comm.~in Probab.}, 4:91--101, 1999.

\bibitem{CSS00}
M.~Cranston, M.~Scheutzow, and D.~Steinsaltz.
\newblock Linear bounds for stochastic dispersion.
\newblock {\em Ann.~Probab.}, 28(4):1852--1869, 2000.

\bibitem{Dim06}
G.~Dimitroff.
\newblock Some properties of isotropic {B}rownian and {O}rnstein-{U}hlenbeck
  flows.
\newblock Ph.D. Dissertation, TU Berlin, \\ URL:
  http://opus.kobv.de/tuberlin/volltexte/2006/1252/, 2006.

\bibitem{DKK04a}
D.~Dolgopyat, V.~Kaloshin, and L.~Koralov.
\newblock A limit shape theorem for periodic stochastic dispersion.
\newblock {\em Comm. Pure Appl. Math.}, 57(9):1127--1158, 2004.

\bibitem{IW81}
N.~Ikeda and S.~Watanabe.
\newblock {\em Stochastic differential equations and diffusion processes}.
\newblock North Holland/Kodansha, Amsterdam/Tokyo, 1981.

\bibitem{Ku90}
H.~Kunita.
\newblock {\em Stochastic Flows and Stochastic Differential Equations}.
\newblock Cambridge University Press, Cambridge, UK, 1990.

\bibitem{LeJ85}
Y.~Le~Jan.
\newblock On isotropic {B}rownian motions.
\newblock {\em Z. Wahrscheinlichkeitstheor. Verw. Geb.}, 70:609--620, 1985.

\bibitem{DaLe88}
Y.~Le~Jan and R.~W.~R. Darling.
\newblock The statistical equilibrium of an isotropic stochastic flow with
  negative {L}yapunov exponents is trivial.
\newblock {\em Lect. Notes Math.}, 1321:175--185, 1988.

\bibitem{MoSch99}
S.-E.~A. Mohammed and M.~K. Scheutzow.
\newblock The stable manifold theorem for stochastic differential equations.
\newblock {\em Ann. Probab.}, 27(2):615--652, 1999.

\bibitem{Oks03}
B.~{\O}ksendal.
\newblock {\em Stochastic differential equations}.
\newblock Springer, Berlin, 6th edition, 2003.

\bibitem{SS02}
M.~Scheutzow and D.~Steinsaltz.
\newblock Chasing balls through martingale fields.
\newblock {\em Ann.~Probab.}, 30(4):2046--2080, 2002.

\bibitem{Wat44}
G.~Watson.
\newblock {\em A treatise on the theory of Bessel functions}.
\newblock Cambridge University Press, Cambridge, 2nd edition, 1944.

\bibitem{Yag57}
A.~M. Yaglom.
\newblock Some classes of random fields in $ n $-dimensional space, related to
  stationary random processes.
\newblock {\em Theory of Probability and its Applications}, 28:273--320, 1957.

\end{thebibliography}

\end{document}